\documentclass[11pt]{article}
\usepackage{amsfonts}
\usepackage{color}
\usepackage{graphics}

\usepackage{cite}
\usepackage{latexsym}
\usepackage{amsmath}
\usepackage{amssymb}
\usepackage[dvips]{epsfig}
\usepackage{amscd}

 \usepackage{color}
\newtheorem{theorem}{Theorem}[section]
\newtheorem{remark}{Remark}[section]

\newtheorem{lemma}[theorem]{Lemma}

\newcommand{\ti}{\tilde}

\def\pf{{\it Proof.}  }

\newcommand{\thatsall}{\hfill$\Box$}
\newcommand{\bi}{\bibitem}

\newcommand{\bt}{\begin{theorem}}
\newcommand{\bl}{\begin{lemma}}
\newcommand{\el}{\end{lemma}}
\newcommand{\et}{\end{theorem}}

\renewcommand{\b}{\beta  }

\newcommand{\te}{\theta}

\newcommand{\al}{\alpha}
\newcommand{\de}{\delta}

\newcommand{\la}{\label}
\newcommand{\si}{\sigma}
\newcommand{\ka}{\kappa}

\newcommand{\bn}{\begin{eqnarray}}
\newcommand{\en}{\end{eqnarray}}
\newcommand{\bnn}{\begin{eqnarray*}}
\newcommand{\enn}{\end{eqnarray*}}

\newcommand{\bnnn}{\begin{eqnarray*}}
\newcommand{\ennn}{\end{eqnarray*}}
\newcommand{\ben}{\begin{enumerate}}
\newcommand{\een}{\end{enumerate}}

\newcommand{\ba}{\begin{aligned}}
\newcommand{\ea}{\end{aligned}}
\newcommand{\be}{\begin{equation}}
\newcommand{\ee}{\end{equation}}

\def\norm[#1]#2{\|#2\|_{#1}}

\def\xiT{\int_0^T}

\makeatletter
\@addtoreset{equation}{section}
\makeatother

\makeatletter
\@addtoreset{equation}{section}
\makeatother

\title{ Nonlinearly Exponential Stability  of
Compressible Navier-Stokes System  with
Degenerate Heat-Conductivity \thanks{ Partially supported by NNSFC   11671027 and 11471321.}
}
\author{Bin Huang, Xiaoding Shi \thanks{
  Email addresses: abinhuang@gmail.com(B. Huang),   shixd@mail.buct.edu.cn (X. Shi)}  \\[3mm]     Department of Mathematics, Faculty  of Science, \\Beijing University of Chemical Technology, \\ Beijing  100029, P. R. China }
\date{ }

\begin{document}
\maketitle

\begin{abstract} We  study the large-time behavior of strong solutions   to the one-dimensional,
compressible Navier-Stokes system for a viscous and heat conducting ideal polytropic gas,
when the viscosity is constant and the heat conductivity
is proportional to a positive power of  the temperature.
Both the specific volume and the temperature are  proved to be  bounded  from below and above independently of time. Moreover, it is shown that the global solution is
nonlinearly exponentially stable   as time tends to infinity.
Note that the conditions imposed on  the  initial data are the same as those of the constant  heat conductivity case ([Kazhikhov-Shelukhin. J. Appl. Math. Mech. 41
(1977); Kazhikhov. Boundary  Value  Problems  for
Hydrodynamical Equations,  50
(1981)] and can be arbitrarily large.
Therefore, our result can be regarded as a natural generalization of the Kazhikhov's ones for the constant heat conductivity case to the degenerate and nonlinear one. \end{abstract}


Math Subject Classification: 35Q35; 76N10.

\section{Introduction}
We consider the compressible Navier-Stokes system,   describing  the  one-dimensional motion of a
viscous heat-conducting    gas, written
  in the Lagrange variables (see \cite{ba,se})\be \la{1.1}
v_t=u_{x},
\ee
\be \la{1.2}
u_{t}+P_{x}=\left(\mu\frac{u_{x}}{v}\right)_{x},\ee
\be  \la{1.'3} \left(e+\frac{u^2}{2}\right)_{t}+ (P
u)_{x}=\left(\kappa\frac{\theta_{x}}{v}+\mu\frac{uu_x}{v}\right)_{x}
 ,
\ee
where   $t>0$ is time,  $x\in (0,1)$ denotes the
Lagrange mass coordinate,  and the unknown functions $v>0,u,$ $\theta>0,e>0,$ and $P$ are,  respectively, the specific volume of the gas,  fluid velocity,   internal energy,  absolute temperature, and  pressure. 
In this paper, we
concentrate on ideal polytropic gas, that is, $P$ and $e$ satisfy \be    P =R \theta/{v},\quad e=c_v\theta +\mbox{ const.},
\ee
where  both specific gas constant  $R$ and   heat
capacity at constant volume $c_v $ are   positive constants. For $\mu$ and $\ka,$
we  consider   the case where $\mu$ and $\ka$ are proportional to (possibly different) powers of  $\te:$  \be   \la{1.3'}\mu=\ti\mu \te^\al, \quad \ka=\ti\ka \te^\beta, \ee
with constants $\ti\mu,\ti\ka>0$ and $\al,\beta\ge 0.$

The system \eqref{1.1}-\eqref{1.3'} is supplemented
with the initial   conditions:
\be(v,u,\te)(x,t=0)=(v_0,u_0,\te_0)(x),  \quad  x\in (0,1), \ee
and boundary   ones:\be \la{1.2a}u(0,t)=u(1,t)=0, \quad  \te_x(0,t)=\te_x(1,t)=0,\quad   t\ge 0.\ee

One can deduce from the Chapman-Enskog expansion for the first level of approximation in kinetic theory that the viscosity $\mu$
and heat conductivity $\ka$ are functions of temperature alone (\cite{cip,cc}). In particular, if the intermolecular
potential varies as $r^{-a},$  with intermolecular distance  $r$, then $\mu$ and $\ka$ are both
proportional to the power $(a+4)/(2a)$ of the temperature, that is, \eqref{1.3'} holds with $\al=\beta=(a+4)/(2a).$ Indeed, for Maxwellian molecules $(a=4),$ the dependence is linear, while for elastic spheres $(a\rightarrow\infty),$ the dependence is like $\te^{1/2}.$

For constant coefficients $(\al=\beta=0)$  and large initial data, Kazhikhov and Shelukhin
\cite{ks} first obtained    the global existence
     of solutions  in bounded domains. From then on,
significant progress has been made on the mathematical aspect of the
initial
boundary value problems, see \cite{az1,az2,az3,akm,ji4,ji3,kaw} and the references therein.
Moreover, much effort has been made to generalize this approach to other cases and in
particular to models satisfying \eqref{1.3'}, which in fact has proved to be challenging especially for temperature dependence on $\mu.$ Motivated by the fact that in the case of isentropic flow a temperature dependence in the viscosity translates
into a density dependence, there is a body of literature (see \cite{hs1,hsj1,akm,jiq1,kaq1,pan1}  and the references therein) studying  the case that $\mu$  is
independent of $\te$, and heat conductivity is allowed to depend on temperature in a special
way with a positive lower bound and balanced with corresponding constitution relations.

When it comes to the
physical case \eqref{1.3'} with $\al=\beta$,  there is few results partially because of   the possible degeneracy and strong nonlinearity in viscosity and heat diffusion  introduced in such relations. As
a first step in this direction, Jenssen-Karper \cite{jk1}  proved the global existence of a weak solution to \eqref{1.1}--\eqref{1.2a}
under the assumption that $\al=0$ and $\beta\in(0,3/2).$ Later, for $\al=0$ and $\beta\in(0,\infty),$  Pan-Zhang \cite{pz1} obtain the global strong solution under the condition that \be \la{1.8}( v_0,u_0,\te_0)\in H^1\times H^2\times H^2. \ee

Concerning  the large-time behavior of the strong solutions to \eqref{1.1}-\eqref{1.2a},    Kazhikhov \cite{kaz1} (see also \cite{akm,az1,az2,az3,na1,na2,na3,ni,qy} among others) first obtains that for the case that $\al=\b=0,$   the strong solution   is
nonlinearly exponentially stable   as time tends to infinity. However, it should be mentioned here that the methods used there relies heavily on the non-degeneracy of the heat conductivity $\ka$ and  cannot be applied directly to  the degenerate and nonlinear case  ($\al=0,\b>0$).  In fact,  one of the main aims of this paper is to show that for $\al=0 $ and $\b>0,$ the global strong solutions obtained by \cite{pz1}  are indeed asymptotically stable as time tends to infinity. Moreover, we will improve  the results of \cite{pz1} by relaxing their assumptions on the initial data. Then  we state
  our main result   as follows.
 \begin{theorem} \la{thm11}Suppose that \be\la{1.9}\al=0,\quad \beta> 0,\ee
 and that the initial data $ ( v_0,u_0,\te_0)$   satisfies
  \be \la{1.10} ( v_0,\te_0)\in   H^1 (0,1),\quad u_0\in H^1_0 (0,1),\ee  and \be\la{1.11}
\inf_{x\in (0,1)}v_0(x)>0, \quad \inf_{x\in (0,1)}\theta_0(x)>0. \ee
Then, the initial-boundary-value problem \eqref{1.1}-\eqref{1.2a} has a unique strong solution $(v,u,\te)$ satisfying\be
  \la{3k}\begin{cases}   v , \,\theta \in L^\infty(0,\infty;H^1(0,1)),\quad u  \in L^\infty(0,\infty;H^1_0(0,1)),\\ v_t\in
  L^\infty(0,\infty;L^2(0,1))\cap L^2(0,\infty;H^1(0,1)), \\ u_x,\,\te_x,\,u_t,\,\theta_t,\,v_{xt},\,u_{xx},\,\theta_{xx} \,\in
  L^2((0,1)\times(0,\infty)).\end{cases}\ee
  Moreover, there exists some positive constants $C$ and $\eta_0>0$ such that  for any $(x,t)\in (0,1)\times(0,\infty),$  \be C^{-1}\le v(x,t)\le C,\quad C^{-1}\le\te(x,t)\le C,\ee and that for any $t>0,$ \be \left\|\left(v-\int_0^1v_0dx,u,\te-\int_0^1\te_0dx\right)(\cdot,t)\right\|_{H^1(0,1)}\le C e^{-\eta_0t}.\ee
  \end{theorem}

A few remarks are in order.
  \begin{remark} For  $\al=\b=0,$ under the conditions   \eqref{1.10} and \eqref{1.11}, Kazhikhov and Shelukhin \cite{ks} first obtained existence   of global strong solutions to the initial-boundary-value problem \eqref{1.1}-\eqref{1.2a}. Later,    Kazhikhov   \cite{kaz1} further proves that   the strong solution   is
nonlinearly exponentially stable   as time tends to infinity.  Therefore, our Theorem \ref{thm11} can be regarded as a natural generalization of the classical results   \cite{ks,kaz1} to the degenerate and nonlinear case that $\al=0,\b>0.$ \end{remark}

  \begin{remark}As far as the existence   of global strong solutions is concerned, our result also improves  Pan and Zhang's result \cite{pz1} where they need the initial data satisfy  \eqref{1.8}   which are stronger than \eqref{1.10}.  \end{remark}

We now make some comments on the analysis of this paper. The key step to study the large-time behavior  of  the global strong  solutions is to get the     time-independent lower and upper bounds   of both $v$ and $\theta $ (see   \eqref{cqq1} and \eqref{tq2}).  Compared with \cite{kaz1,ni}, the main difficulty comes from the degeneracy and nonlinearity of the heat conductivity because of  $\b>0.$      Hence, to obtain \eqref{tq2}, some new ideas are needed. The key observations are as follows: First, after using the standard   energetic  estimate (see \eqref{etvt}) and modifying  the idea due to Kazhikhov \cite{kaz1}, we obtain that the specific volume $v$ is bounded from above and below time-independently (see \eqref{cqq1}).  Then, although it seems difficult to obtain the uniform  lower bound of $\te$  at first, after   observing that (see \eqref{qw2s})  \bnn \xiT\max\limits_{x\in[0,1]}\left| \te^{1/2}(x,t)- \int_0^1\te^{1/2}(x,t)dx\right|^2 dt\le C ,\enn   we prove  that the $L^\infty(0,\infty;L^p)$-norm of $\te^{-1}$ is bounded (see \eqref{ljj1}),  which in turn not only implies that   (see \eqref{lm1} and \eqref{lm2})\bnn \xiT\max\limits_{x\in[0,1]}\left| \te (x,t)- \int_0^1\te (x,t)dx\right|^2 dt+\xiT\int_0^1 u_x^2dxdt\le C,\enn but also yields that
  the $L^2((0,1)\times (0,T))$-norm of $\te_x$ is bounded provided $\b>1$ (see \eqref{sq1}). Finally, for $\b\in(0,1], $ we find that  the $L^2((0,1)\times (0,T))$-norm of $\te_x$ can be  bounded by  the $L^4( 0,T;L^2(0,1))$-norm of $u_x$ which plays an important role in obtaining  the uniform bound on $L^2((0,1)\times (0,T))$-norm of both $\te_x$ and $u_{xx} $ (see Lemma \ref{lemm5}) for $\b\in (0,1]$.   The whole procedure will be carried out in the next section.
\section{ Proof of Theorem \ref{thm11}}

We first state   the following the local existence result which can be proved by using the principle of compressed mappings (c.f. \cite{kan,nash,tan}).

\begin{lemma} \la{lma1} Let \eqref{1.9}-\eqref{1.11} hold. Then there exists some $T>0$ such that  the initial-boundary-value problem \eqref{1.1}-\eqref{1.2a} has a unique strong solution $(v,u,\te)$ satisfying\bnn
 \begin{cases}   v , \,\theta  \in L^\infty(0,T;H^1(0,1)), \quad u  \in L^\infty(0,T;H^1_0(0,1)), \\ v_t\in
  L^\infty(0,T;L^2(0,1))\cap  L^2(0,T;H^1(0,1)), \\ u_t,\,\theta_t,\,v_{xt},\,u_{xx},\,\theta_{xx} \,\in
  L^2((0,1)\times(0,T)).\end{cases}\enn \end{lemma}

Then,   the   a priori   estimates (see \eqref{cqq1},  \eqref{vx1}, \eqref{lm5}, \eqref{tq2}, and \eqref{tq7} below) where the constants  depend only on the data of the problem   make it possible to continue the local solution      to the whole interval $[0,\infty)$ and finish the proof of Theorem \ref{thm11}.

Next, without loss of generality, we assume that $\ti\mu= \ti\ka=R=c_v=1 $ and that
 \be \la{1.3a} \int_0^1 v_0dx=1, \quad \int_0^1 \left(\frac{u_0^2}{2}+\te_0\right)dx=1.\ee
  Motivated by the second
law of thermodynamics, one has   the following standard energetic  estimate   embodying the dissipative effects of viscosity and thermal diffusion.
\begin{lemma} \la{2k} It holds that\be\ba\la{etvt} \sup_{0\le t<\infty}\int_0^1\left(  \frac{u^2}{2}+ (v-\ln
v)+ (\theta-\ln \theta)\right)dx+\int_0^t V(s)ds \le E_0,  \ea\ee
where
\bnn\ba V(t)\triangleq\int_0^1\left(\frac{\theta^\b\theta_{x}^2}{v\te^2}+\frac{u_{x}^{2}}{v\te}
\right)(x,t)dx, \ea\enn and \bnn E_0\triangleq\int_0^1\left(\frac{u_0^2}{2}+ (v_0-\ln v_0)+ (\theta_0-\ln
\theta_0)\right)dx.\enn

\end{lemma}

\pf It follows from  \eqref{1.1}, \eqref{1.'3}, \eqref{1.2a},  and   \eqref{1.3a}  that for $t>0$
 \be\la{wq1}  \int_0^1 v(x,t)dx=  1,  \quad\int_0^1 \left(\frac{u^2}{2}+\te\right)(x,t)dx= 1.\ee

 Noticing that the energy equation \eqref{1.'3} can be written as
\be\la{tee1}\theta_{t}+  \frac\theta v
u_{x}= \left(\frac{\theta^\b\theta_{x}}{v}\right)_{x}+ \frac{u_{x}^{2}}{v},\ee    multiplying   \eqref{1.1}   by
 $ 1- {v}^{-1} $,  \eqref{1.2}   by $
u$,   (\ref{tee1}) by
$  1- {\theta}^{-1} $,  and adding them altogether, we get
\bnn  \ba&\left(u^2/2+ (v-\ln v)+ (\theta-\ln
\theta)\right)_t+ \frac{u^2_x}{v\theta}
+ \frac{\te^\b\theta_x^2}{v\theta^2}\\&=\left(\frac{  u u_x}{v}-\frac{  u
\theta}{v}\right)_x + u_x+ \left(\left(1-\theta^{-1}\right)
\frac{\te^\b\theta_x}{v}\right)_x,\ea\enn  which together with \eqref{1.2a} yields \eqref{etvt} and finishes the proof of Lemma \ref{2k}. \thatsall

Next, we derive the following representation of $v$ which is essential in obtaining the time-independent upper and lower bounds of $v$.

\begin{lemma} \la{lma1} We have the following expression of $v$
\be \la{qw9}\ba v(x,t)=   D(x,t) Y(t) + \int_0^t\frac{D(x,t) Y(t)\te(x,\tau)}{ D(x,\tau) Y(\tau)}d\tau , \ea\ee where \be \la{qw10}\ba D(x,t)=&v_0(x)\exp\left\{\int_{0}^x \left(u(y,t)-u_0(y)\right)dy\right\}\\&\times \exp\left\{-\int_0^1v\int_0^xudydx+\int_0^1v_0\int_0^xu_0dydx\right\}, \ea\ee and
\be\la{qw11}\ba Y(t)=\exp\left\{-\int_0^t\int_0^1\left(u^2+\te\right) dxds\right\} .\ea\ee
\end{lemma}

\pf First, denoting by
\be\la{wqq1}\ba  \si\triangleq\frac{u_x}{v}-\frac{\te}{v} , \ea\ee
it follows from  \eqref{1.2} that
\be\la{yq1}\ba \left(\int_0^xudy\right)_t=\si-\si(0,t) , \ea\ee
 which implies
\bnn\ba v\si(0,t)=v\si-v\left(\int_0^xudy\right)_t . \ea\enn
Integrating this in $x $ over $(0,1)$ together with \eqref{wq1} and \eqref{wqq1} yields\be\la{yq2}\ba    \si(0,t) &=\int_0^1v\si dx-\int_0^1v\left(\int_0^xudy\right)_tdx\\&= \int_0^1\left(u_x-\te\right) dx-\left(\int_0^1v\int_0^xudydx\right)_t+\int_0^1u_x\int_0^xudy dx\\&=- \left(\int_0^1v\int_0^xudydx\right)_t-\int_0^1\left(\te+ u^2 \right) dx.\ea\ee

Next, since $u_x=v_t,$
 we have \bnn \si=(\ln v)_t-\frac{\te}{v},\enn which together with \eqref{yq1} and \eqref{yq2} gives
\bnn\ba \left(\int_0^xudy\right)_t  = (\ln v)_t-\frac{\te}{v}+ \left(\int_0^1v\int_0^xudydx\right)_t+\int_0^1\left(\te+ u^2 \right) dx . \ea\enn
Integrating this over $(0,t)$ leads to
\be\la{qw8}\ba v(x,t) = D(x,t) Y(t) \exp\left\{\int_0^t\frac{\te}{v}ds\right\}, \ea\ee
 with $D(x,t)$ and $Y(t)$ as in  \eqref{qw10} and \eqref{qw11}  respectively.

Finally, denoting by
\bnn\ba  g =\int_0^t\frac{\te}{v}ds,\ea\enn
we have by using \eqref{qw8}
\bnn\ba g_t=\frac{\te(x,t)}{v(x,t)}=\frac{\te(x,t)}{ D(x,t) Y(t)\exp \{g\}} ,\ea\enn
which gives
\bnn\ba  \exp\{g\}  =1+\int_0^t\frac{\te(x,
\tau)}{ D(x,\tau) Y(\tau)}d\tau .\ea\enn Putting this into \eqref{qw8} leads to \eqref{qw9}.
\thatsall

With Lemmas \ref{lma1}  and \ref{2k} at hand,   we are in a position to prove  the time-independent upper and lower bounds of $v$.
\begin{lemma}\la{wq20} For any $(x,t)\in[0,1]\times [0,+\infty),$ it holds
\be\ba\la{cqq1}
C^{-1}\le v(x,t)\le C,
\ea\ee where (and in what follows)   $C $  denotes some
generic positive constant
 depending only on $\b,\|(v_0-1,u_0,\theta_0-1)\|_{H^1(0,1)},
 \inf\limits_{x\in [0,1]}v_0(x),$ and $ \inf\limits_{x\in [0,1]}\theta_0(x).$
\end{lemma}
\pf First, since the function $x-\ln x $ is convex, Jensen's inequality gives
\bnn\ba \int_{0}^{1}  \te dx-\ln\int_{0}^{1}   \theta dx  \le \int_{0}^{1} (\te-\ln\theta ) dx ,\ea\enn
  which together with \eqref{etvt} and \eqref{wq1} leads to \be \la{qwe1} \bar \te(t)\triangleq \int_0^{1}\te(x,t)dx\in[ \al_1,1],\ee   where $0<\al_1<\al_2$  are two roots of  \bnn x-\ln x =E_0.\enn

  Next, both \eqref{wq1} and Cauchy's inequality imply
  \bnn\ba
 \left|\int_0^1 v\int_0^x udydx\right|&\le \int_0^1 v\left|\int_0^x udy\right|dx\\&\le\int_0^1 v\left(\int_0^1 u^2dy\right)^{1/2}dx\\&\le C,
 \ea\enn which combined with \eqref{qw10} gives
\be\ba \la{cq8}C^{-1}\le  D(x,t)\le  C.\ea\ee

Furthermore,  one deduces from \eqref{wq1} that\bnn\ba\la{cq9} 1\le \int_0^1\left(u^2+\te\right) dx\le 2,\ea\enn
which   yields that for any $0\le \tau<t<\infty, $
\be \la{cq7}e^{-2t}\le Y(t)\le 1,\quad e^{-2(t-\tau)}\le  \frac{Y(t)}{Y(\tau)}\le e^{- (t-\tau)}.\ee

Next, it follows from \eqref{wq1}  that
\bnn\ba\la{bqq2}  \left|\te^{\frac{\b+1}{2}}(x,t)- \bar\te^{\frac{\b+1}{2}}( t) \right| &\le \frac{\b+1}{2}\left(\int_0^1 \frac{ \te^\b \te_x^2}{\te^2 v} dx\right)^{1/2}\left(\int_0^1  {\te v} dx\right)^{1/2}\\ &\le C V^{1/2}(t)  \max_{x\in [0,1]}v^{1/2}(x,t),\ea\enn which together with \eqref{qwe1} leads to
\be\ba\la{bqq7}
 \frac{\alpha_1}{4}-CV(t)\max_{x\in [0,1]}v(x,t)\le  \theta(x,t)\le C +CV(t)\max_{x\in [0,1]}v(x,t),
\ea\ee
for all $(x,t)\in [0,1]\times [0,\infty).$

Next, it follows from \eqref{qw9}, \eqref{cq8},  \eqref{cq7}, and \eqref{bqq7} that
\bnn\ba\la{cq10} v(x,t)
&\le C+C\int_0^t e^{-t+\tau}\max_{x\in[0,1]} \te (x,\tau)d\tau\\
&\le C+C\int_0^t e^{-t+\tau} \left(1+V(\tau)\max_{x\in [0,1]}v(x,\tau)\right)d\tau\\
&\le C+C\int_0^t  V(\tau)\max_{x\in [0,1]}v(x,\tau) d\tau,\ea\enn
which together with  the Gronwall inequality gives
\be\ba\la{cq11}  v(x,t) \le C \ea\ee  for all $(x,t)\in[0,1]\times[0,+\infty).$ Combining this with \eqref{qw9}, \eqref{cq8},  \eqref{cq7}, and \eqref{bqq7} yields that
\be\ba\la{cq12a} v(x,t)&\ge C \int_0^t e^{-2(t-\tau)}\min_{x\in[0,1]}\te(x,\tau)d\tau \\&\ge C \int_0^t e^{-2(t-\tau)}\left(\frac{\al_1}{4}-C  V( \tau)\right)d\tau\\ &\ge \frac{C \al_1}{8} - \frac{C \al_1}{8} e^{-2t} -C   \int_0^t e^{-2(t-\tau)} V( \tau) d\tau.\ea\ee
Noticing that
\bnn\ba  \int_0^t e^{-2(t-\tau)} V( \tau) d\tau&=\int_0^{t/2} e^{-2(t-\tau)} V( \tau) d\tau+\int_{t/2}^t e^{-2(t-\tau)} V( \tau) d\tau\\&\le e^{-t}\int_0^\infty V(\tau)d\tau +\int_{t/2}^t V( \tau) d\tau\rightarrow 0, \mbox{ as }t\rightarrow \infty,\ea\enn we deduce from \eqref{cq12a} that  there exists some $\ti T >0$ such that
\be\ba\la{cq12} v(x,t) \ge \frac{C \al_1}{16} \ea\ee for all $(x,t)\in[0,1]\times[ \ti T ,+\infty).$

 Finally,    using \eqref{qw8}, \eqref{cq8}, and \eqref{cq7},  we obtain that  there exists some positive constant $C $ such that \bnn\ba v(x,t) \ge C^{-1}\ea\enn for all $(x,t)\in[0,1]\times[0, \ti T ].$
Combining this, \eqref{cq12}, and \eqref{cq11}   gives  \eqref{cqq1} and finishes the proof of Lemma \ref{wq20}.
\thatsall

To obtain the uniform (with respect to time) lower bound of the temperature, we  need the following   time-independent bound on the $L^\infty(0,T;L^p)$-norm of $\te^{-1}.$

\begin{lemma} \la{len1} For any $p>0,$   there exists some positive constant $C(p)$  such that \be \la{ljj1}  \sup_{0\le t\le T}\int_0^1 \te^{1-p}dx+\int_0^T\int_0^1 \frac{\te^\b\te_x^2}{\te^{p+1}}dxdt+\int_0^T\int_0^1\frac{u_x^2}{\te^p}dxdt  \le C(p).\ee\end{lemma}
 \pf First, it follows from \eqref{etvt} that  \eqref{ljj1} holds for $p=1.$

 Next, for $p\not=1,$ multiplying \eqref{tee1}   by $1/\te^p$    and integration by parts gives
\be\ba\la{aa1}  &\frac{1}{p-1}\left(\int_0^1\te^{1-p}dx\right)_t+p\int_0^1 \frac{\te^\b\te_x^2}{v\te^{p+1}}dx+\int_0^1\frac{u_x^2}{v\te^p}dx \\& =\int_0^1 \frac{\left(\te^{1-p}-1\right)u_x}{v}dx+\int_0^1 \frac{ u_x}{v}dx\\ & \le C(p)\int_0^1   \left|\te^{\frac{1}{2}}-1 \right|\left(1+\te^{\frac{1}{2}-p}\right)| u_x| dx +\left(\int_0^1 \ln v dx\right)_t  \\ & \le  C(p)\max_{x\in[0,1]}\left| \te^{  \frac{1}{2}}-1\right|\left( \int_0^1 |u_x|dx+\left( \int_0^1 \te^{1-p}dx \right)^{1/2}\left( \int_0^1\frac{u_x^2}{v\te^p}dx \right)^{1/2} \right)  \\&\quad  +\left(\int_0^1 \ln v dx\right)_t\\ & \le  C(p)\max_{x\in[0,1]}\left| \te^{  \frac{1}{2}}-1\right|^2+C(p) \left( \int_0^1 |u_x|dx\right)^2 + \frac12 \int_0^1\frac{u_x^2}{v\te^p}dx \\&\quad  + C(p)\max_{x\in[0,1]}\left| \te^{  \frac{1}{2}}-1\right|^2\int_0^1 \te^{1-p}dx  +\left(\int_0^1 \ln v dx\right)_t.\ea\ee

 Next, we claim that for any real number $q,$ there exists a positive constant $C(q)$ such that
\be\la{te1}
1-\bar{\te}^{q}\le C (q)V^{1/2}(t) .  \ee

Indeed, standard calculation  gives
\be\ba\la{24eq1}  1-\bar \te^{q} &=\int_0^1 \frac{d}{d\eta }\left(\left(\int_0^1(\te+\eta  \frac{u^2}{2})dx\right)^{q}\right) d\eta  \\ &=q\int_0^1 \left(\int_0^1 (\te+\eta  \frac{u^2}{2})dx\right)^{ q-1} d\eta \cdot \int_0^1\frac{u^2}{2} dx \\&\le C(q) \int_0^1 {u^2}  dx ,\ea\ee where in the last inequality we have used  $$\alpha_1\le\int_0^1 \te dx\le \int_0^1(\te+\eta\frac{u^2}{2})dx\le 1,$$ due to   \eqref{wq1}  and \eqref{qwe1}.
Furthermore, \be\ba\la{24eq2}   \int_0^1 {u^2}  dx\le  \max_{x\in[0,1]} |u|\left(\int_0^1 {u^2}  dx\right)^{1/2} \le C \int_0^1 |u_x|dx,  \ea\ee and
\be \la{ux}\ba  \int_0^1 |u_x|dx   &\le  \left( \int_0^1\frac{ u_x^2}{v\te}dx\right)^{1/2}\left(\int_0^1 v\te dx\right)^{1/2}\\ &\le C\left( \int_0^1\frac{ u_x^2}{v\te}dx\right)^{1/2}\le CV^{1/2}(t). \ea\ee
Combining all these estimates \eqref{24eq1}--\eqref{ux} gives \eqref{te1}.

Then,   it follows from \eqref{te1} and \eqref{qwe1} that for $ \b\in (0,1),$
\be\la{qwq2}\ba \max_{x\in[0,1]}\left|\te^{\frac{1 }{2}}-1\right|&\le \max_{x\in[0,1]}\left|\te^{\frac{1 }{2}}-\bar\te^{\frac{1 }{2}}\right| +\max_{x\in[0,1]}\left|\bar\te^{\frac{1 }{2}}-1\right| \\ &\le  C  \int_0^1  \te^{ -\frac{1 }{2} }|\te_x|dx+C V^{1/2}(t)\\&\le C\left(\int_0^1  \te^{\b-2}\te_x^2dx\right)^{1/2}\left(\int_0^1  \te^{1-\b} dx\right)^{1/2}+ CV^{1/2}(t)\\&\le CV^{1/2}(t),\ea\ee
and that for $\b\ge 1,$
 \be\la{qw2}\ba \max_{x\in[0,1]}\left|\te^{\frac{1 }{2}}-1\right|&\le \max_{x\in[0,1]}\left|\te^{\frac{1 }{2}}-\bar\te^{\frac{1 }{2}}\right| +\max_{x\in[0,1]}\left|\bar\te^{\frac{1 }{2}}-1\right|\\ &\le C \max_{x\in[0,1]}\left|\te^{\frac{\b }{2}}-\bar\te^{\frac{\b}{2}}\right| +CV^{1/2}(t) \\ &\le  C  \int_0^1  \te^{ \frac{\b }{2}-1}|\te_x|dx+C V^{1/2}(t)\\&\le  CV^{1/2}(t). \ea\ee
It thus follows from \eqref{qwq2}, \eqref{qw2}, and \eqref{etvt} that
 \be\la{qw2s}\ba\int_0^T \max_{x\in[0,1]}\left|\te^{\frac{1 }{2}}-1\right|^2dt\le C.\ea\ee

Finally, noticing  that  for $p\in(0,1),$
\bnn\la{qwq1}\ba
\int_0^1\te^{1-p}dx\le \int_0^1\te dx+ 1\le C,
\ea\enn and that both
  \eqref{etvt} and   \eqref{wq1} imply
\be\la{22eq5}
\sup_{0\le t<\infty}\int_0^1  |\ln v|dx\le C,
\ee
 after using \eqref{ux}, \eqref{etvt}, \eqref{qw2s},   and Gronwall's inequality,   we   obtain  \eqref{ljj1} from \eqref{aa1} and   finish  the proof of Lemma \ref{len1}.
\thatsall

Next, using Lemma \ref{len1}, we have the following estimate on the $L^\infty(0,T;L^2)$-norm of $v_x.$

\begin{lemma}\la{lemm3} There exists a   positive constant $C $  such that
\be\la{vx1}\ba  \sup_{0\le t\le T}\int_0^1 v_x^2dx+\int_0^T\int_0^1 v_x^2(\te+1) dxdt\le C, \ea\ee
for any $T\ge 0.$  \end{lemma}

\pf
First, choosing $p=\b$ in \eqref{ljj1} gives
\be  \la{lmm3}\int_0^T  \int_0^1  \te^{-1 }   \te_x^2dx  dt  \le C,\ee which together with \eqref{qwe1} implies
\be\ba\la{lm1}  \int_0^T \max_{x\in[0,1]}\left(\te (x,t) -\bar\te ( t)\right)^2dt    \le C\int_0^T  \int_0^1  \te^{-1 }   \te_x^2dx\int_0^1  \te    dx dt  \le C.  \ea\ee

Next, integrating the momentum equation  \eqref{1.2}  multiplied by $u$ over $[0,1] $ with respect to $x,$  we obtain after  integrating by parts
\be\la{bg1} \ba & \frac12\left(\int_0^1 u^2dx\right)_t+\int_0^1 \frac{u_x^2}{v}dx \\&=\int_0^1 \frac{\te }{v}u_xdx \\&=\int_0^1 \frac{\left(\te-\bar\te\right)}{v}u_xdx+\left(\bar\te -1\right)\int_0^1 \frac{ u_x}{v}dx + \int_0^1 \frac{ u_x}{v}dx \\&\le C \max_{x\in[0,1]}\left(\te-\bar\te\right)^2+CV(t) + \left(\int_0^1 \ln vdx\right)_t, \ea\ee where in the last inequality we have used \eqref{te1} and \eqref{ux}. Combining this with  \eqref{etvt}, \eqref{22eq5}, and \eqref{lm1} yields \be\ba\la{lm2} \int_0^T\int_0^1 u_x^2dxdt\le C. \ea\ee

Next, using  \eqref{1.1}, we
   rewrite momentum equation  \eqref{1.2} as
 \be\ba\la{mom1}
\left(\frac{v_x}{v}\right)_t=u_t+\left(\frac{\te}{v}\right)_x,
\ea\ee due to
\bnn
\left(\frac{v_t}{v}\right)_x=\left(\frac{v_x}{v}\right)_t .
\enn
Multiplying   \eqref{mom1}  by $\frac{v_x}{v}$ leads to
 \be\ba\la{mom2}
\frac{1}{2}\left[\left(\frac{v_x}{v}\right)^2\right]_t&=\frac{v_x}{v}u_t+\frac{v_x}{v}\left(\frac{\te}{v}\right)_x
\\&=\left(\frac{v_x}{v}u\right)_t-u(\ln v)_{xt}+\frac{v_x\te_x}{v^2}-\frac{v_x^2\te}{v^3}
\\&=\left(\frac{v_x}{v}u\right)_t-\left[u(\ln v)_{t}\right]_x+\frac{u_x^2}{v}+\frac{v_x\te_x}{v^2}-\frac{v_x^2\te}{v^3}.\ea\ee
Integrating \eqref{mom2} over $[0,1]\times[0,T],$   one has
 \bnn\ba\la{mom3}
&\sup_{0\le t \le T}\int_0^1\left[\frac{1}{2}\left(\frac{v_x}{v}\right)^2-\frac{v_x}{v}u\right]dx
+\int_0^T\int_0^1\frac{v_x^2\te}{v^3}dxdt\\
&\le C+\int_0^T\int_0^1 \frac{v_x\te_x}{v^2}dxdt\\
&\le C+\frac12\int_0^T\int_0^1\frac{v_x^2\te}{v^3}dxdt+ C\int_0^T\int_0^1  \te^{-1} \te_x^2dxdt\\
&\le C+\frac12\int_0^T\int_0^1\frac{v_x^2\te}{v^3}dxdt,
\ea\enn due to \eqref{lmm3}. This in particular implies
 \be\la{vvx1}\ba
\sup_{0\le t \le T}\int_0^1 v_x^2  dx+\int_0^T\int_0^1 v_x^2\te dxdt\le C,\ea\ee due to the following simply fact: \bnn\la{mom4}\ba
\int_0^1\frac{v_x}{v}udx\le \frac{1}{4}\int_0^1\left(\frac{v_x}{v}\right)^2dx+ C.
\ea\enn

Finally, it follows from \eqref{vvx1} that
\bnn\ba  \bar\te \int_0^1 v_x^2dx&= \int_0^1 v_x^2 \left(\bar\te-\te\right) dx+  \int_0^1 v_x^2\te dx\\ &\le \frac{\bar \te}{2}\int_0^1 v_x^2  dx +\frac{1}{2\bar \te}\int_0^1 v_x^2  \left(\te-\bar\te\right)^2 dx+  \int_0^1 v_x^2\te dx\\& \le \frac{\bar \te}{2}\int_0^1 v_x^2  dx +C\max_{x\in[0,1]}\left(\te-\bar\te\right)^2 +  \int_0^1 v_x^2\te dx,\ea\enn which  together with \eqref{lm1} and \eqref{vvx1} leads to
\bnn \int_0^T\int_0^1 v_x^2  dxdt\le C.\enn
 Combining  this with \eqref{vvx1} gives \eqref{vx1} and finishes the proof of
 Lemma \ref{lemm3}. \thatsall

For further uses, we need the following estimate on the $L^2((0,1)\times(0,T))$-norm of $\te_x$ for $\b\in (0,1].$

 \begin{lemma}\la{lemm6} If $ 0<\b\le 1,$ there exists a   positive constant $C $ such that
  \be\ba\la{sq2} \int_0^T\int_0^1  \te_x^2dxdt
\le C+C \int_0^T\left(\int_0^1 u_x^2dx\right)^2dt ,\ea\ee for any $T>0.$\end{lemma}
\pf    Multiplying \eqref{tee1}   by $\te^{1-\frac{\b}{2}}$ and integration by parts gives
\be \la{p2}\ba  &\frac{2}{4-\b}\left(\int_0^1\te^{2-\frac{\b}{2}}dx\right)_t +\frac{(2-\b)}{2}\int_0^1 \frac{\te^{\frac{\b}{2}}\te_x^2}{v}dx\\& =-\int_0^1 \frac{\te^{2-\frac{\b}{2}} }{v}u_xdx+\int_0^1\frac{\te^{1-\frac{\b}{2}} u_x^2}{v} dx \\&=\int_0^1 \frac{\left(\bar\te^{2-\frac{\b}{2}}-\te^{2-\frac{\b}{2}}\right) }{v}u_xdx + \left(1-\bar\te^{2-\frac{\b}{2}} \right)\int_0^1 \frac{ u_x}{v}dx\\&\quad-\int_0^1 \frac{ u_x}{v}dx+\int_0^1\frac{\te^{1-\frac{\b}{2}} u_x^2}{v} dx\\ &\le C\int_0^1  \left|\te^{2-\frac{\b}{2}}-\bar\te^{2-\frac{\b}{2} } \right|\left| u_x\right|dx  +CV(t)-\left(\int_0^1 \ln v dx\right)_t+\int_0^1\frac{\te^{1-\frac{\b}{2}} u_x^2}{v} dx,\ea\ee where in the last inequality we have used \eqref{te1}.
Direct calculation yields that
\be\la{p4}\ba &\int_0^1  \left|\te^{2-\frac{\b}{2}}-\bar\te^{2-\frac{\b}{2} } \right|\left| u_x\right|dx \\&\le C\max_{x\in[0,1]}\left|\te^{1-\frac{\b}{4}}-\bar\te^{1-\frac{\b}{4}} \right|\left(\int_0^1\left(\te^{2-\frac{\b}{2}}+1\right)dx \right)^{1/2}\left(\int_0^1u_x^2dx\right)^{1/2}\\&\le C\left(\int_0^1\te^{-\frac{\b}{4}}|\te_x|dx\right)^2+C \int_0^1\left(\te^{2-\frac{\b}{2}}+1\right)dx \int_0^1u_x^2dx \\&\le C \int_0^1\te^{ 1-\frac{\b}{2}}  dx  \int_0^1\te^{-1} \te_x^2dx +C \int_0^1\left(\te^{2-\frac{\b}{2}}+1\right)dx \int_0^1u_x^2dx\\&\le C   \int_0^1\te^{-1} \te_x^2dx +C \int_0^1\left(\te^{2-\frac{\b}{2}}+1\right)dx \int_0^1u_x^2dx ,\ea\ee
 and that for any $\delta>0$\be\la{p5}\ba&\int_0^1\frac{\te^{1-\frac{\b}{2}} u_x^2}{v} dx \\ &\le C\left(\max_{x\in[0,1]}\left|\te^{ 1-\frac{\b}{2} }-\bar \te^{1-\frac{\b}{2} }\right|+1\right)   \int_0^1    u_x^2  dx  \\&\le C\int_0^1 \te^{-\frac{\b}{2}}|\te_x|dx\int_0^1   u_x^2  dx+C\int_0^1   u_x^2  dx  \\&\le \delta\int_0^1\left(\te^{ -1}+\te^{\frac{\b}{2}}\right) \te_x^2dx+C(\delta)\left(\int_0^1   u_x^2  dx\right)^2+C\int_0^1   u_x^2 dx. \ea\ee
Putting \eqref{p4}  and \eqref{p5}  into \eqref{p2},   choosing $\de$ suitably small, and using \eqref{lm2}, \eqref{lmm3}, and the Gronwall inequality, one obtains \bnn\la{p1}\ba   \int_0^1\te^{2-\b/2}dx  +\int_0^T\int_0^1 \te^{\b/2}\te_x^2 dxdt\le C+C \int_0^T\left(\int_0^1 u_x^2dx\right)^2dt,\ea\enn which together with \eqref{etvt} implies
  \bnn\ba  \int_0^T\int_0^1  \te_x^2dxdt& \le  \int_0^T\int_0^1  \left(\te^{\b-2}+\te^{\b/2}\right) \te_x^2dxdt\\&\le C+\int_0^T\int_0^1  \te^{\b/2} \te_x^2dxdt\\&
\le C+C \int_0^T\left(\int_0^1 u_x^2dx\right)^2dt .\ea\enn
 This gives \eqref{sq2} and finishes the proof of Lemma \ref{lemm6}.\thatsall

Then, we have the following uniform estimate on the $L^2((0,1)\times (0,T))$-norm of $u_t$ and $u_{xx}.$

\begin{lemma}\la{lemm5} There exists a   positive constant $C $ such that
\be\ba\la{lm5}  \sup_{0\le t\le T}\int_0^1 u_x^2dx+\int_0^T\int_0^1(u_t^2+ u_{xx}^2)dxdt\le C.\ea\ee for any $T\ge 0.$\end{lemma}
\pf First, we rewrite the momentum equation  \eqref{1.2}  as
\be\ba\la{mom5} u_t-\frac{u_{xx}}{v}=-\frac{u_xv_x}{v^2}-\frac{\te_x}{v } +\frac{\te v_x}{v^2}.\ea\ee
Multiplying both sides of \eqref{mom5} by $u_{xx}$ and integrating the resultant equality in $x$ over $[0,1]$ lead to
\be\la{qu2}\ba
&\frac{1}{2}\frac{d}{dt}\int_0^1 u_x^2dx+\int_0^1\frac{u_{xx}^2}{v} dx\\
&\le \left|\int_0^1\frac{u_xv_x}{v^2}u_{xx}dx\right|+\left|\int_0^1\frac{\te_x}{v }u_{xx}dx\right| +\left|\int_0^1\frac{\te v_x}{v^2}u_{xx}dx\right|\\
&\le \frac{1}{4}\int_0^1\frac{u_{xx}^2}{v} dx+C\int_0^1\left(u_x^2v_x^2+v_x^2\te^2 +\te_x^2\right)dx.
\ea\ee
Direct computation yields that for any $\de>0,$
\be\la{qu1}\ba &\int_0^1 \left(u_x^2v_x^2+v_x^2\te^2+\te_x^2\right)dx\\ &\le C\left(\max_{x\in[0,1]} u_x^2+\max_{x\in[0,1]}\left(\te-\bar\te\right)^2+1\right)\int_0^1  v_x^2dx +\int_0^1 \te_x^2dx\\ &\le C\max_{x\in[0,1]} u_x^2+C\max_{x\in[0,1]}\left(\te-\bar\te\right)^2  +C\int_0^1  v_x^2dx +\int_0^1 \te_x^2dx\\ &\le \delta \int_0^1 u_{xx}^2dx+ C(\delta) \int_0^1 u_x^2dx+C\max_{x\in[0,1]}\left(\te-\bar\te\right)^2  +C\int_0^1  v_x^2dx +\int_0^1 \te_x^2dx,\ea\ee
where in the last inequality  we have used
\be\ba\la{ux2}
\max_{x\in[0,1]}u_x^2&\le \int_0^1 \left|\left(u_x^2\right)_x\right|dx\\&\le 2\left(\int_0^1 u_{xx}^2d x\right)^{1/2}\left(\int_0^1 u_x^2dx\right)^{1/2} \\&\le \delta \int_0^1 u_{xx}^2dx+ C(\delta) \int_0^1 u_x^2dx,
\ea\ee due to   $\int_0^1 u_xdx=0.$
Putting \eqref{qu1} into \eqref{qu2} and choosing $\de$ suitably small yields
\be\ba\la{le5eq1}  \int_0^1 u_x^2dx +\int_0^T\int_0^1 u_{xx}^2 dxdt\le C   +C\int_0^T \int_0^1\te_x^2dx dt,\ea\ee due to \eqref{lm1}, \eqref{lm2}, and \eqref{vx1}.

Next, on the one hand, if $\b>1,$ choosing $p=\b-1$ in \eqref{ljj1} gives
\be\ba\la{sq1} \int_0^T\int_0^1  \te_x^2dxdt \le C, \ea\ee
 which  along with    \eqref{le5eq1}  gives
\be\ba\la{lm5w}  \sup_{0\le t\le T}\int_0^1 u_x^2dx+\int_0^T\int_0^1 u_{xx}^2dxdt+\int_0^T\int_0^1\te_x^2dxdt\le C.\ea\ee
 On the other hand, if $\b\in (0,1],$ it follows from \eqref{le5eq1}, \eqref{sq2},  \eqref{lm2},   and Gronwall's inequality that \eqref{lm5w}  still holds.

Finally, it follows from \eqref{mom5}, \eqref{lm5w}, \eqref{qu1}, \eqref{lm1}, \eqref{lm2}, and \eqref{vx1} that
\bnn  \int_0^T\int_0^1 u_t^2dxdt\le C,\enn  which together with \eqref{lm5w} gives \eqref{lm5} and finishes the proof of Lemma \ref{lemm5}.\thatsall

Now, we can prove the uniform lower and upper bounds of the temperature $\te.$

\begin{lemma}\la{lmm8} There exists a positive constant $C$ such that for any $(x,t)\in[0,1]\times[0,T]$ \be\ba\la{tq2}  C^{-1}\le \te(x,t)\le C  . \ea\ee\end{lemma}

\pf First, for $p>\b+1,$
multiplying \eqref{tee1} by
$ \te^{p-1}$ and integrating the resultant equality in $x$ over $(0,1)$  leads to
 \be\la{p7}\ba&\frac{1}{p}\left(\int_0^1 \te^pdx\right)_t  +(p-1)\int_0^1 \frac{\te^{p+\b-2}\te_x^2}{v}dx \\& =\int_0^1 \frac{\te^{p-1}u_x^2}{v}dx-\int_0^1 \frac{\te^{p}u_x }{v}dx\\&\le C\max_{x\in[0,1]} u_x^2\int_0^1  \te^{p-1} dx+C\int_0^1 |(\te^{p}-1)u_x|dx -\int_0^1 \frac{  u_x }{v}dx\\& \le C\max_{x\in[0,1]} u_x^2\int_0^1  \te^{p-1} dx+C\int_0^1 |\te-1|\left(\te^{p-1}+1\right)|u_x|dx -\int_0^1 \frac{  u_x }{v}dx  \\& \le C\max_{x\in[0,1]} \left(u_x^2+|\te-1|^2\right)\left(1+\int_0^1  \te^{p} dx\right)  -\left(\int_0^1 \ln vdx\right)_t.\ea\ee

It follows from \eqref{lm5w}, \eqref{ux2},   and  \eqref{lm2} that
\be \la{uux2}\int_0^T \max_{x\in[0,1]}  u_x^2dt+\int_0^T\int_0^1 \te_x^2 dx dt\le C,\ee which together with \eqref{te1} shows \be \la{p8}\ba\int_0^T\max_{x\in[0,1]}  |\te-1|^2dt&\le C\int_0^T\left(\max_{x\in[0,1]}  \left|\te-\bar\te\right|^2+\max_{x\in[0,1]}  \left|\bar\te-1\right|^2\right)dt\\&\le C\int_0^T\int_0^1 \te_x^2dxdt+C\int_0^TV(t)dt\le C.\ea \ee
 Combining \eqref{p7}-\eqref{p8} with the Gronwall inequality   gives \be\ba\la{lm7eq1}\sup_{0\le t\le T}\int_0^1 \te^pdx+\int_0^T\int_0^1 \te^{p+\b-2}\te_x^2dxdt\le C(p) .\ea\ee

 Next, multiplying \eqref{tee1} by
$  \te^{\b}\te_t$ and integrating the resultant equality over $(0,1)$ yields
\bnn\ba &  \int_0^1 \te^\b\te_t^2dx+ \int_0^1\frac{ \te^{\b+1}\theta_tu_x}{v  }dx\\&=\int_0^1\te^{\b}\te_t\left(\frac{\te^\b\theta_{x}}{v}\right)_{x}dx+\int_0^1\frac{ \te^{\b }\theta_tu_x^2}{v }dx  \\& =-\int_0^1\frac{\te^\b\theta_{x}}{v}\left(\te^{\b}\te_t\right)_{x}dx+\int_0^1\frac{ \te^{\b }\theta_tu_x^2}{v }dx   \\& =-\int_0^1\frac{\te^\b\theta_{x}}{v}\left(\te^{\b}\te_x\right)_{t}dx+\int_0^1\frac{ \te^{\b }\theta_tu_x^2}{v }dx  \\&=-\frac{1}{2}
\int_0^1\frac{\left((\te^\b\theta_{x})^2\right)_t}{v}dx+\int_0^1\frac{ \te^{\b }\theta_tu_x^2}{v }dx \\&=-\frac{1}{2} \left(\int_0^1\frac{(\te^\b\theta_{x})^2}{v}dx\right)_t
-\frac{1}{2}\int_0^1\frac{(\te^\b\theta_{x})^2u_x}{v^2}dx+\int_0^1\frac{ \te^{\b }\theta_tu_x^2}{v }dx ,\ea\enn
which gives
\be\ba\la{lm8eq2} &  \int_0^1 \te^\b\te_t^2dx+ \frac{1}{2} \left(\int_0^1\frac{(\te^\b\theta_{x})^2}{v}dx\right)_t\\&=- \frac{1}{2}\int_0^1\frac{(\te^\b\theta_{x})^2u_x}{v^2}dx-\int_0^1\frac{ \te^{\b+1}\theta_t u_x}{v  }dx+\int_0^1\frac{ \te^{\b }\theta_tu_x^2}{v }dx\\&\le C\max_{x\in[0,1]} |u_x|\int_0^1\left(\te^\b\te_x\right)^2dx +\frac{1}{2}\int_0^1 \te^\b\te_t^2dx +C\int_0^1 \te^{\b+2}u_x^2dx\\&\quad+C\int_0^1\te^\b u_x^4dx\\&\le C\left(\int_0^1\left(\te^\b\te_x\right)^2dx\right)^2+C \max_{x\in[0,1]} u_x^2+C \max_{x\in[0,1]} u_x^4+\frac{1}{2}\int_0^1 \te^\b\te_t^2dx  \ea\ee
due to \eqref{lm7eq1}.

Next, it follows from  \eqref{ux2}  and \eqref{lm5}    that
\be\la{tq3} \int_0^T \max_{x\in[0,1]}  u_x^4dt\le C,\ee which together with \eqref{lm8eq2},  the Gronwall inequality, \eqref{uux2}, and \eqref{lm7eq1} leads to
\be\ba\la{tebtex}  \sup_{0 \le t\le T}\int_0^1 \left(\te^\b\theta_{x}\right)^2 dx+\int_0^T\int_0^1 \te^\b\te_t^2dxdt\le C. \ea\ee
This in particular gives
\bnn\ba \max_{x\in[0,1]}\left|\te^{\b+1}-\bar\te^{\b+1}\right|& \le (\b+1)\int_0^1 \te^\b|\te_x|dx\\&\le C\left(\int_0^1 \left(\te^\b\te_x\right)^2dx\right)^{1/2} \le C,\ea\enn which implies that for all $(x,t)\in[0,1]\times [0,\infty),$
\be\la{teup} \te(x,t)\le C.\ee

Next, it follows from \eqref{te1} that
\be\ba \la{lm8eq5}&\int_0^T \int_0^1 \left(\te^{\b+2}-1\right)^2dxdt\\ &\le 2\int_0^T \int_0^1 \left(\te^{\b+2}-\bar \te^{\b+2}\right)^2dxdt+C\int_0^T  V(t)dt\\ &\le C\int_0^T \left(\int_0^1  \te^{{\b+1}}|\te_x| dx\right)^2dt+C
\\ &\le C\int_0^T \int_0^1  \te^{\b-2} \te_x^2 dx dt+C  \le C,\ea\ee where in the last inequality we have used \eqref{teup}. Combining this, \eqref{tebtex}, and \eqref{teup} in particular gives
\bnn\ba\la{lm8eq6} &\int_0^T \left|\left(\int_0^1 \left(\te^{\b+2}-1\right)^2dx\right)_t\right|dt\\&=2\int_0^T \left| \int_0^1 \left(\te^{\b+2}-1\right) \left(\te^{\b+2}\right)_tdx\right|dt\\ &\le C\int_0^T  \int_0^1 \left(\te^{\b+2}-1\right)^2dxdt+C\int_0^T\int_0^1   \te^{2\b+2}\te^2_tdx dt\\&\le C +C\int_0^T\int_0^1   \te^{\b}\te^2_tdx dt \le C,\ea\enn which together with
 \eqref{lm8eq5}   leads to
\be\la{kl1}\ba \lim_{t\rightarrow \infty}\int_0^1 \left(\te^{\b+2}-1\right)^2dx=0. \ea\ee

Then, we claim that
\be\la{kl2}\ba  \lim_{t\rightarrow \infty}\left(1-\bar \te\right)=0,\ea\ee
which combined with \eqref{kl1} gives
  \be\ba\la{tetebar} \lim_{t\rightarrow \infty}\int_0^1 \left(\te^{\b+2} -\bar \te^{\b+2}\right)^2dx=0. \ea\ee
It thus follows from  \eqref{tebtex}  and  \eqref{teup}  that
\bnn\ba & \max_{x\in [0,1]}\left(\te^{\b+2}-\bar\te^{\b+2}\right)^2\\ &\le C\int_0^1\left|\te^{\b+2}-\bar\te^{\b+2}\right|\left|(\te^{\b+2})_x\right|dx \\ &\le C\left(\int_0^1\left(\te^{\b+2} -\bar \te^{\b+2}\right)^2dx\right)^{1/2}\left(\int_0^1 \left(\te^\b\te_x\right)^2dx\right)^{1/2}  \\ &\le C\left(\int_0^1\left(\te^{\b+2} -\bar \te^{\b+2}\right)^2dx\right)^{1/2}  ,  \ea\enn
which together with \eqref{tetebar} and \eqref{kl2}  implies that there exists some  $  T_0>0$ such that
\be\la{x5}\ba \te (x,t)\ge 1/2, \ea\ee   for all $(x,t)\in [0,1]\times [ T_0,\infty).$ Moreover, it follows from  \cite[Lemma 2.2]{pz1} that  there exists some constant $C \ge  2$
 such that  \bnn \theta(x,t)\ge  C^{-1} ,\enn  for all $(x,t)\in [0,1]\times [0,T_0].$ Combining this, \eqref{x5}, and \eqref{teup} gives \eqref{tq2}.

Finally,   it   remains to prove \eqref{kl2}. Indeed,  it follows from \eqref{bg1} and \eqref{teup}  that
\bnn \la{ul1}\ba   \frac12\left(\int_0^1 u^2dx\right)_t+\int_0^1 \frac{u_x^2}{v}dx =\int_0^1 \frac{\te}{v}u_xdx\le C\int_0^1 |u_x|dx ,\ea\enn
which yields that there exists some constant $C$ such that for any $N>0 $ and $s,t\in [N,N+1]$
\bnn\ba \int_0^1 u^2(x,t)dx-\int_0^1 u^2(x,s)dx\le C\int_N^{N+1}\int_0^1 |u_x| dxd\tau .\ea\enn
Integrating  this with respect to $s$ over $(N, N+1)$ and using \eqref{24eq2} yields that
\be\la{kll1}\ba \sup_{N\le t\le N+1} \int_0^1 u^2(x,t)dx \le  C\left(\int_N^{N+1}\int_0^1  u^2_x  dxd\tau\right)^{1/2}.\ea\ee
Letting $N\rightarrow \infty$ in \eqref{kll1}   gives
\bnn\ba  \lim_{t\rightarrow \infty}\int_0^1 u^2(x,t)dx=0,\ea\enn
due to \eqref{lm2}. Combining this   with \eqref{wq1} gives \eqref{kl2} and finishes the proof of Lemma  \ref{lmm8}. \thatsall

Next,
 we have the following uniform estimate on the $L^2((0,1)\times (0,T))$-norm of $\te_t$ and $\te_{xx}.$

\begin{lemma}\la{tq8}There exists a positive constant $C$ such that \be\la{tq7}\ba \sup_{0\le t\le T}\int_0^1 \te_x^2dx+\int_0^T\int_0^1\left( \te_t^2+\te_{xx}^2\right)dxdt\le C.\ea\ee\end{lemma}
\pf First, both \eqref{tq2} and \eqref{tebtex} lead  to
\be\ba\la{lm9eq2}  \sup_{0 \le t\le T}\int_0^1 \theta_{x}^2 dx+\int_0^T\int_0^1 \te_t^2dxdt\le C. \ea\ee

Next, it follows from \eqref{tee1} that
\bnn\ba \frac{\left(\te^\b\te_x\right)_x}{v}=  \frac{ \te^\b\te_xv_x}{v^2}- \frac{u_x^2}{v}+ \frac{ \te u_x}{v}+\te_t,\ea\enn
which together with \eqref{vx1}, \eqref{tq2},  \eqref{lm2}, \eqref{tq3},   and \eqref{lm9eq2} gives
\be\la{tq4}\ba  \int_0^T\int_0^1\left| \left(\te^\b\te_x\right)_x\right|^2dxdt  & \le C\int_0^T\max_{x\in[0,1]}\left(  \te^\b\te_x\right)^2\int_0^1 v_x^2dxdt+C\\ & \le C\int_0^T\max_{x\in[0,1]}\left(  \te^\b\te_x\right)^2 dt+C.\ea\ee Since $\te_x(0,t)=0,$ we get by \eqref{uux2}  and \eqref{tq2}, \bnn\ba \int_0^T\max_{x\in[0,1]}\left(  \te^\b\te_x\right)^2 dt &\le  \int_0^T\int_0^1\left| \left(\left(\te^\b\te_x\right)^2\right)_x\right| dxdt \\& \le   C   (\de) +\de \int_0^T\int_0^1\left| \left(\te^\b\te_x\right)_x\right|^2dxdt ,  \ea\enn  which together with  \eqref{tq4} and \eqref{tq2}  implies
\be\la{tq5}\ba
&\int_0^T \max_{x\in[0,1]} \te_x^2dt+\int_0^T\int_0^1\left| \left(\te^\b\te_x\right)_x\right|^2dxdt  \le C.
\ea\ee

Finally, since
\bnn\te_{xx}=\frac{\left(\te^\b\te_x\right)_x}{\te^\b} -  \frac{\b \te_x^2}{\te } ,\enn it follows from \eqref{tq5},  \eqref{tq2},  and \eqref{tebtex} that
\bnn \ba\int_0^T\int_0^1  \te_{xx}^2dxdt  &\le  C    \int_0^T\int_0^1\left| \left(\te^\b\te_x\right)_x\right|^2dxdt+ C \int_0^T\max_x \te_x^2\int_0^1 \te_x^2dxdt  \\&\le C+C\sup_{0\le t\le T}\int_0^1 \te_x^2dx \int_0^T\max_{x\in[0,1]}\te_x^2dt
\\&\le  C, \ea\enn  which together with \eqref{lm9eq2}
gives \eqref{tq7} and finishes the proof of Lemma \ref{tq8}.\thatsall

Finally, we have the following nonlinearly exponential stability of the strong solutions.
\begin{lemma} There exist  some positive constants $C$ and $\eta_0 $ both
 depending only on $\b,\|(v_0-1,u_0,\theta_0-1)\|_{H^1(0,1)},
 \inf\limits_{x\in [0,1]}v_0(x),$ and $ \inf\limits_{x\in [0,1]}\theta_0(x) $ such that
\be \la{klq1} \|(v-1,u,\te-1)(\cdot,t)\|_{H^1(0,1)}\le Ce^{-\eta_0t}.\ee\end{lemma}
\pf Noticing that all the constants $C$ in Lemmas \ref{lemm3}, \ref{lemm5}, and \ref{tq8}    are independent of $T,$ we have
\be\la{zma1} \int_0^\infty \left(\left|\frac{d}{dt}\|v_x(\cdot,t)\|_{L^2}^2\right|+\left|\frac{d}{dt}\|u_x(\cdot,t)\|_{L^2}^2\right| +\left|\frac{d}{dt}\|\te_x(\cdot,t)\|_{L^2}^2\right|\right)dt\le C,\ee where we have used \bnn
 \int_0^1u_xu_{xt}dx=- \int_0^1u_tu_{xx}dx.\enn It thus follows from \eqref{zma1} that
\bnn \lim_{t\rightarrow \infty}\|(v_x,u_x,\te_x)(\cdot,t)\|_{L^2(0,1)}=0,\enn which in particular implies
\bnn \lim_{t\rightarrow \infty}\|(v-1,u ,\te-1)(\cdot,t)\|_{H^1(0,1)}=0.\enn
Therefore, since we know that the temperature remains bounded from above and below independently of time and the solution becomes small in $H^1$-norm for large time $t$,
we can conclude the solution decays to the constant state exponentially as $t\rightarrow \infty,$  that is, \eqref{klq1} holds (c.f. \cite{ok}). \thatsall

 \end{document}